\documentclass[14pt]{article}
\usepackage{mathrsfs}
\usepackage{amsthm}
\usepackage{amssymb}
\usepackage{amsmath}
\usepackage{graphicx}
\usepackage{color}
\usepackage{amsfonts}
\usepackage{float}
\usepackage{cite}
\usepackage[text={140mm,210mm},left=45mm,vmarginratio=1:1]{geometry}
\newtheorem{theorem}{Theorem}[section]

\newtheorem{lemma}[theorem]{Lemma}

\numberwithin{equation}{section}
\normalsize

\begin{document}
\title{\textbf{Asymptotic for critical value of the large-dimensional SIR epidemic on clusters}}

\author{Xiaofeng Xue \thanks{\textbf{E-mail}: xfxue@bjtu.edu.cn \textbf{Address}: School of Science, Beijing Jiaotong University, Beijing 100044, China.}\\ Beijing Jiaotong University}

\date{}
\maketitle

\noindent {\bf Abstract:} In this paper we are concerned with the
SIR (Susceptible-Infective-Removed) epidemic on open clusters of
bond percolation on the squared lattice. For the SIR model, a
susceptible vertex is infected at rate proportional to the number of
infective neighbors while an infective vertex becomes removed at a
constant rate. A removed vertex will never be infected again. We
assume that there is only one infective vertex at $t=0$ and define
the critical value of the model as the maximum of the infection
rates with which infective vertices die out with probability one,
then we show that the critical value is $\big(1+o(1)\big)/(2dp)$ as
$d\rightarrow+\infty$, where $d$ is the dimension of the lattice and
$p$ is the probability that a given edge is open. Our result is a
counterpart of the main theorem in \cite{Hol1981} for the contact
process.

\quad


\section{Introduction}\label{section one}

In this paper, we are concerned with the SIR
(susceptible-infective-removed) epidemic model on open clusters of
bond percolation in squared lattices $\{\mathbb{Z}^d\}_{d\geq 1}$
(see a survey of percolation in \cite{Grimmett1999}). For later use,
we identify $\mathbb{Z}^d$ with the vertices set of it and denote by
$E_d$ the edges set of $\mathbb{Z}^d$. We denote by $O$ the origin
of the lattice. We assume that $\{X(e)\}_{e\in E_d}$ are i. i. d.
random Bernoulli variables such that
\[
P\big(X(e)=1\big)=p=1-P\big(X(e)=0\big)
\]
for some $p\in (0,1]$. For later use, we write $X(e)$ as $X(x,y)$
when $e$ connects vertices $x$ and $y$. For vertices $x$ and $y$, we
write $x\sim y$ when and only when there is an edge $e$ connecting
$x,y$ and $X(e)=1$. Intuitively, we delete each edge in state $0$
while remain those in state $1$, then $x\sim y$ when and only when
they are neighbors on the consequent graph.

We denote by $\mathcal{P}(\mathbb{Z}^d)$ the set of all the subsets
of $\mathbb{Z}^d$, then the SIR model is a Markov process with state
space
\[
\Omega=\Big\{(A,B):~A,B\in \mathcal{P}(\mathbb{Z}^d), A\bigcap
B=\emptyset\Big\}.
\]
We denote by $(S_t, ~I_t)$ the state of the process at moment $t$
for any $t\geq 0$, then the SIR epidemic evolves as follows.
\begin{equation}\label{equ 1.1}
(S_t,~I_t)\rightarrow
\begin{cases}
&(S_t,~I_t\setminus\{x\}) \text{\quad at rate~}1 \text{~if~}x\in
I_t,\\
&(S_t\setminus\{x\},~I_t\bigcup \{x\})\text{\quad at
rate~}\lambda\sum\limits_{y:y\sim x}\mathbf{1}_{\{y\in
I_t\}}\text{~if~} x\in S_t,
\end{cases}
\end{equation}
where $\lambda$ is a positive parameter called the infection rate
and we denote by $\mathbf{1}_A$ the indicator function of the random
event $A$.

Intuitively, the process $\{(S_t, ~I_t)\}_{t\geq 0}$ describes the
spread of an epidemic. Vertices in $S_t$ are susceptible which can
be infected while vertices in $I_t$ are infective which can infect
neighbors. Vertices in $\mathbb{Z}^d\setminus \big(S_t\bigcup
I_t\big)$ are removed which will never be infected again. A
susceptible vertex is infected at rate proportional to the number of
infective neighbors while an infective vertex becomes removed at
rate one. Note that here we say $x$ and $y$ are neighbors when the
edge $e$ connecting them satisfies $X(e)=1$ as we introduced.

The main topic we are concerned with in this paper is the estimation
of the critical value of our model, which is the maximum of the
infection rates with which infective vertices die out with
probability one when at $t=0$ there are finite infective vertices.
The critical value of infection is first studied for another type of
epidemic which is the SIS model, where an infective vertex will
become healthy and then may be infected again. The SIS model is also
named as the contact process. See a survey of the contact process in
Chapter 6 of \cite{Lig1985} and Part one of \cite{Lig1999}. A direct
corollary of our main result given in the next section can be seen
as a counterpart of the asymptotic behavior of critical value of the
large-dimensional contact process obtained in \cite{Hol1981} by
Holley and Liggett. For mathematical details, see the next section.

Let $p_c$ be the maximum of $p$ with which the open cluster
containing $O$ is finite with probability one, then when $p<p_c$
infective vertices die out almost surely since the infection spreads
on finite graphs and the critical value of infection rate is
infinity as a result. However, Kesten proves
$\lim\limits_{d\rightarrow+\infty}2dp_c(d)=1$ in \cite{Kes1990} and
hence $p>p_c(d)$ for given $p>0$ and sufficiently large $d$, which
makes the critical infection rate of the epidemic nontrivial in
large dimension. We are inspired by the technique introduced in
\cite{Kes1990} a lot when proving the main result of this paper.

We are inspired a lot by recent references about the SIR epidemic on
percolation models. The percolation on complete graph is also known
as the ER (Erd\H{o}s-R\'{e}nyi) graph. In \cite{Neal2003}, Neal
studies a discrete-time version of SIR on the ER graph and gives
limit distribution of the process. In \cite{Xue2016a}, Xue considers
a law of large numbers of the SIR on ER graph inspired by the theory
of density dependent population model introduced by Ethier and Kurtz
in \cite{Ethier1986}. In \cite{Xue2015} and \cite{Xue2016b}, Xue
considers the SIR epidemics on open clusters of oriented site and
oriented bond percolation models on lattices as auxiliary tools to
study corresponding contact processes.

\section{Main result}\label{section two}
In this section we give our main result. First we introduce some
notations, definitions and basic assumptions. For each $d\geq 1$, we
assume that $\{X(e)\}_{e\in E_d}$ are defined under the probability
space $\big(Y_d,\mathcal{F}_d,\mu_d\big)$. We denote by $E_{\mu_d}$
the expectation operator with respect to $\mu_d$. For any $\omega\in
Y_d$, we denote by $P_{\lambda}^\omega$ the probability measure of
the process $\{(S_t,~I_t)\}_{t\geq 0}$ with infection rate $\lambda$
in the random environment on $\mathbb{Z}^d$ with respect to
$\{X(\omega,~e)\}_{e\in E_d}$. $P_{\lambda}^\omega$ is called the
quenched measure. We denote by $E_\lambda^\omega$ the expectation
operator with respect to $P_{\lambda}^\omega$. We define
\[
P_{\lambda,d}\big(\cdot\big)=E_{\mu_d}\Big(P_{\lambda}^{\omega}\big(\cdot\big)\Big)=\int P_\lambda^\omega(\cdot)~\mu_d(d\omega),
\]
which is called the annealed measure. We denote by $E_{\lambda,d}$
the expectation operator with respect to $P_{\lambda,d}$. When there
is no misunderstanding, we write $Y_d, \mathcal{F}_d, \mu_d,
E_{\mu_d}, P_{\lambda,d}, E_{\lambda,d}$ as $Y,\mathcal{F},\mu, ,
E_\mu, P_\lambda, E_\lambda$.

Throughout this paper we assume that
\begin{equation}\label{equ 2.1 basic assumption}
\big(S_0,~I_0\big)=\big(\mathbb{Z}^d\setminus \{O\},~\{O\}\big)
\end{equation}
for the process on $\mathbb{Z}^d$ and each $d\geq 1$. Note that $O$
is the origin of the lattice as we introduced. According to the
basic coupling of Markov processes (see Section 3.1 of
\cite{Lig1985}), for any $\lambda_1<\lambda_2$,
\[
P_{\lambda_1}\big(I_t\neq \emptyset,\forall ~t>0\big)\leq
P_{\lambda_2}\big(I_t\neq \emptyset,\forall ~t>0\big).
\]
As a result, it is reasonable to define
\begin{equation}\label{equ 2.2}
\lambda_c(d)=\sup\big\{\lambda:~P_{\lambda,d}\big(I_t\neq
\emptyset,\forall~t>0\big)=0\big\}
\end{equation}
for each $d\geq 1$. That is to say, $\lambda_c$ is the maximum of
the infection rates with which the infective vertices die out almost
surely when there are finite infective vertices at $t=0$.

The following theorem is our main result, which gives the asymptotic
behavior of $\lambda_c(d)$ as $d$ grows to infinity.

\begin{theorem}\label{theorem 2.1 main}
If $\lambda_c(d)$ is defined as in Equation \eqref{equ 2.2}, then
\[
\lim_{d\rightarrow+\infty}d\lambda_c(d)=\frac{1}{2p}.
\]
\end{theorem}

Note that $p$ is the probability that a given edge is in state $1$
as we defined at the beginning of this paper. When $p=1$, our model
reduces to the classic SIR epidemic on lattices and Theorem
\ref{theorem 2.1 main} shows that
\[
\lim_{d\rightarrow+\infty}2d\lambda_c(d)=1.
\]
Let $\widehat{\lambda}_c(d)$ be the counterpart of $\lambda_c(d)$
with respect to the contact process, then it is proved in
\cite{Hol1981} that
\[
\lim_{d\rightarrow+\infty}2d\widehat{\lambda}_c(d)=1
\]
for the classic case where $p=1$. For general case, it is easy to
see that $\widehat{\lambda}_c(d)\leq \lambda_c(d)$ according to
basic coupling of Markov processes. Hence as a direct corollary of
Theorem \ref{theorem 2.1 main},
\[
\limsup_{d\rightarrow+\infty}d\widehat{\lambda}_c(d)\leq
\frac{1}{2p}.
\]
This result has been proved in \cite{Xue2016c} in a general case
where each infective vertex recovers at i. i. d. random rates. We
believe that
$\liminf\limits_{d\rightarrow+\infty}d\widehat{\lambda}_c(d)\geq\frac{1}{2p}$
and hence
$\lim\limits_{d\rightarrow\infty}d\widehat{\lambda}_c(d)=\frac{1}{2p}$
but have not found a proof yet.

We give an intuitive explanation of Theorem \ref{theorem 2.1 main}
according to a mean-field analysis. When $d$ is large, an vertex has
about $2dp$ neighbors according to the law of large numbers. Each
infective vertex becomes removed at rate one while infects a given
neighbor at rate $\lambda$, hence the number $I$ of infective
vertices approximately follows the ODE
\[
\frac{dI}{dt}=(2dp\lambda-1)I.
\]
Then $I$ converges to $0$ when and only when
$\lambda<\frac{1}{2dp}$.

The proof of Theorem \ref{theorem 2.1 main} is divided into two
sections. In Section \ref{section three}, we will show that
\begin{equation}\label{equ 2.3}
\limsup_{d\rightarrow+\infty}d\lambda_c(d)\leq \frac{1}{2p}.
\end{equation}
Our proof of Equation \eqref{equ 2.3} is inspired by the technique
introduced in \cite{Kes1990}. We consider the self-avoiding paths on
which the infection spreads from the beginning to the end. We call
such paths the infection paths. If $\lambda$ satisfies that the
probability that there exist infection paths with arbitrary lengths
is positive, then $\lambda$ is a upper bound of $\lambda_c$. In
Section \ref{section four}, we will show that
\begin{equation}\label{equ 2.4}
\liminf_{d\rightarrow+\infty}d\lambda_c(d)\geq \frac{1}{2p}.
\end{equation}
The proof of Equation \eqref{equ 2.4} is relative easy. We consider
the number $N$ of the vertices which have ever been infective. If
$\lambda$ satisfies that the mean of $N$ is finite, then $\lambda$
is a lower bound of $\lambda_c$. For mathematical details, see
Sections \ref{section three} and \ref{section four}.

\section{Proof of Equation \eqref{equ 2.3}}\label{section three}
In this section we will give the proof of Equation \eqref{equ 2.3}. Since Equation \eqref{equ 2.3} is about the asymptotic behavior of $\lambda_c(d)$ as $d\rightarrow+\infty$,  we assume that the dimension $d$ of the lattice satisfies $d\geq 20$ throughout this section. We will explain the reason of this assumption later. First
we introduce some definitions and notations. We denote by
$\|\cdot\|_1$ the $l_1$ norm on $\mathbb{Z}^d$ such that
\[
\|x-y\|_1=\sum_{i=1}^d|x_i-y_i|
\]
for any $x=(x_1,\ldots,x_d),~y=(y_1,\ldots,y_d)\in \mathbb{Z}^d$.
For $1\leq i\leq d$, we use $e_i$ to denote
\[
(0,\ldots,0,\mathop 1\limits_{i \text{th}},0,\ldots,0),
\]
which is the $i$th elementary vector on $\mathbb{Z}^d$.
For any $x\in Z^d$, we assume that $T(x)$ is an exponential time
with rate $1$. For any $x,y\in \mathbb{Z}^d$ such that
$\|x-y\|_1=1$, we assume that $U(x,y)$ is an exponential time with
rate $\lambda$. Note that here we care about the order of $x$ and
$y$, hence $U(x,y)\neq U(y,x)$. We assume that all these exponential
times are independent and are independent with the random
environment $\{X(e)\}_{e\in E_d}$. Intuitively, $T(x)$ is the time
interval $x$ waits for to become removed after being infected while
$U(x,y)$ is the time interval $x$ waits for to infect $y$ after $x$
being infected if $X(x,y)=1$.

For each integer $K\geq 1$, we define
\begin{align*}
L_K=\Big\{\vec{l}&=(l_0,l_1,\ldots,l_K)\in
\big(\mathbb{Z}^d\big)^K:~l_0=O;~\|l_{i+1}-l_i\|_1=1,\\
&\forall~0\leq i\leq K-1;~l_i\neq l_j,\forall~i\neq j\Big\}
\end{align*}
as the set of self-avoiding paths on $\mathbb{Z}^d$ starting at the
origin $O$ with length $K$. For each $K\geq 1$ and any
$\vec{l}=(l_0,\ldots,l_K)\in L_K$, we define
\[
A_{\vec{l}}=\big\{U(l_{i},~l_{i+1})\leq T(l_i),\forall~0\leq i\leq
K-1;~X(l_i,l_{i+1})=1, \forall~0\leq i\leq K-1\big\}.
\] 
Note that $A_{\vec{l}}$ is a random event. According to the
definition of the SIR model, a vertex can not be infected repeatedly
and hence any vertex that has ever been infective must be infected
through a self-avoiding path from $O$ to it since $O$ is the only
infective vertex at $t=0$. Therefore, in the sense of coupling,
\begin{equation}\label{equ 3.1}
A_{\vec{l}}\subseteq \big\{l_K\in I_t \text{~for some~}t>0\big\}
\end{equation}
for $\vec{l}=(l_0,\ldots,l_K)\in L_K$ and
\begin{equation}\label{equ 3.2}
\big\{x\in I_t\text{~for
some~}t>0\big\}=\bigcup\limits_{K=1}^{+\infty}\bigcup_{\vec{l}\in
L_K,\atop l_K=x}A_{\vec{l}}
\end{equation}
for any $x\neq O$. By direct calculation, for any $\vec{l}\in L_K$,
\begin{equation}\label{equ 3.3}
P_\lambda(A_{\vec{l}})=\Big(P(U\leq T)P\big(X(e)=1\big)\Big)^K=\frac{\lambda^Kp^K}{(\lambda+1)^K},
\end{equation}
where $U$ and $T$ are independent exponential times with rates $\lambda$ and $1$ respectively, since all the edges on the self-avoiding path are different with each other and $\{X(e)\}_{e\in E_d}$ are i. i. d.. For later use, we need to give an upper bound of $P(A_{\vec{l}}\bigcap A_{\vec{s}})$ for $\vec{l},\vec{s}\in L_K$. For this purpose, for $\vec{l}=(l_0,\ldots,l_K), \vec{s}=(s_0,\ldots,s_K)\in L_K$, we define
\[
D(\vec{l},~\vec{s})=\Big\{0\leq i\leq K:~s_i=l_j\text{~for some~}j\in \{0,\ldots,K\}\Big\}
\]
and
\[
F(\vec{l},~\vec{s})=\Big\{0\leq i\leq K-1:~s_i=l_j \text{~and~}s_{i+1}=l_{j+1}\text{~for some~}j\in \{0,\ldots,K\}\Big\}.
\]
We use ${\rm card}(A)$ or $|A|$ to denote the cardinality of the set $A$, then $|D(\vec{l},~\vec{s})|$ is the number of vertices that both $\vec{l}$ and $\vec{s}$ visit while $|F(\vec{l},~\vec{s})|$ is the number of edges that $\vec{l}$ and $\vec{s}$ visit through the same direction. We have the following lemma which gives an upper bound of $P(A_{\vec{l}}\bigcap A_{\vec{s}})$.
\begin{lemma}\label{lemma 3.1}
For $\vec{l},\vec{s}\in L_K$,
\[
P(A_{\vec{l}}\bigcap A_{\vec{s}})\leq \big(\frac{\lambda p}{\lambda+1}\big)^{2K-|F(\vec{l},~\vec{s})|}\big(\frac{2}{p}\big)^{|D(\vec{l},~\vec{s})\setminus F(\vec{l},~\vec{s})|}~.
\]
\end{lemma}

\proof

By Equation \eqref{equ 3.3},
\[
P(A_{\vec{l}}\bigcap A_{\vec{s}})=P\big(A_{\vec{s}}\big|A_{\vec{l}}\big)P(A_{\vec{l}})=\big(\frac{\lambda p}{\lambda+1}\big)^KP\big(A_{\vec{s}}\big|A_{\vec{l}}\big).
\]
For any $i\not \in D(\vec{l},~\vec{s})$, $X(s_i,s_{i+1})$, $T(s_i)$ and $U(s_i,~s_{i+1})$ are independent with $A_{\vec{l}}$, therefore $P\big(A_{\vec{s}}\big|A_{\vec{l}}\big)$ has the factor
\[
\Big(P(U\leq T)P\big(X(e)=1\big)\Big)^{K-|D(\vec{l},~\vec{s})|}=\big(\frac{\lambda p}{\lambda+1}\big)^{K-|D(\vec{l},~\vec{s})|}.
\]
For each $i\in D(\vec{l},~\vec{s})\setminus F(\vec{l},~\vec{s})$, there exist $0\leq j\leq K-1$ and $x,y,z$ such that $s_i=l_j=x$, $s_{i+1}=y$, $l_{j+1}=z$ and $y\neq z$. Hence there is a factor at most
\[
P\big(U(x,~y)<T(x)\big|U(x,~z)<T(x)\big)\leq \frac{2\lambda}{\lambda+1}
\]
in the expression of $P\big(A_{\vec{s}}\big|A_{\vec{l}}\big)$ for each $i\in D(\vec{l},~\vec{s})\setminus F(\vec{l},~\vec{s})$. Therefore, $P\big(A_{\vec{s}}\big|A_{\vec{l}}\big)$ has a factor at most $\big(\frac{2\lambda}{\lambda+1}\big)^{|D(\vec{l},~\vec{s})\setminus F(\vec{l},~\vec{s})|}$. For each $i\in F(\vec{l},~\vec{s})$, $X(s_i,~s_{i+1})=1$ and $U(s_i,~s_{i+1})\leq T(s_i)$ occurs with probability one conditioned on $A_{\vec{l}}$ since there exists $j$ such that $s_i=l_j$ and $s_{i+1}=l_{j+1}$. In conclusion,
\begin{align*}
P\big(A_{\vec{s}}\big|A_{\vec{l}}\big)&\leq \big(\frac{\lambda p}{\lambda+1}\big)^{K-|D(\vec{l},~\vec{s})|}\big(\frac{2\lambda}{\lambda+1}\big)^{|D(\vec{l},~\vec{s})\setminus F(\vec{l},~\vec{s})|}\\
&=\big(\frac{\lambda p}{\lambda+1}\big)^{K-|F(\vec{l},~\vec{s})|}\big(\frac{2}{p}\big)^{|D(\vec{l},~\vec{s})\setminus F(\vec{l},~\vec{s})|}
\end{align*}
and hence
\[
P(A_{\vec{l}}\bigcap A_{\vec{s}})=\big(\frac{\lambda p}{\lambda+1}\big)^KP\big(A_{\vec{s}}\big|A_{\vec{l}}\big)\leq \big(\frac{\lambda p}{\lambda+1}\big)^{2K-|F(\vec{l},~\vec{s})|}\big(\frac{2}{p}\big)^{|D(\vec{l},~\vec{s})\setminus F(\vec{l},~\vec{s})|}.
\]

\qed

Inspired by the approach introduced in \cite{Kes1990} by Kesten, we consider a special type of self-aoviding paths on $\mathbb{Z}^d$. For each $k\geq 1$,
we define
\begin{align*}
\Gamma_k=\Big\{\vec{l}&=(l_0,\ldots,l_{k\lfloor\log d\rfloor})\in L_{k\lfloor\log d\rfloor}:l_{i+1}-l_i\in \{\pm e_j:1\leq j\leq d-\lfloor\frac{d}{\log d}\rfloor\}\\
&\text{~for any~}i\text{~such that~}\lfloor\log d\rfloor \nmid(i+1); l_{i+1}-l_i\in \{e_j:d-\lfloor \frac{d}{\log d}\rfloor+1\leq j\leq d\}\\
&\text{~for any~}i\text{~such that~}\lfloor\log d\rfloor \mid(i+1)\Big\},
\end{align*}
where we use $a\mid b$ to denote that $b$ is divisible by $a$ and $\{e_j\}_{1\leq j\leq d}$ are the elementary vectors on $\mathbb{Z}^d$ as we defined in Section \ref{section one}. We introduce a random walk $\{S_n\}_{n=0}^{+\infty}$ on $\mathbb{Z}^d$ with paths in $\bigcup_{k\geq 1}\Gamma_k$. For $n=0$, we assume that $S_0=O$. For $n\geq1$, $\{S_n\}_{n=1}^{+\infty}$ evolves as follows. For each $j\geq 1$ such that $\lfloor\log d\rfloor\nmid j$, assuming that we have already obtained the first $j$ steps $S_0, S_1, S_2, \ldots, S_{j-1}$, then
\[
P\big(S_j=z\big|S_i,0\leq i\leq j-1\big)=\frac{1}{|H(j)|}
\]
for any $z\in H(j)$, where
\[
H(j)=\Big\{y: y-S_{j-1}\in \{\pm e_l:1\leq l\leq d-\lfloor\frac{d}{\log d}\rfloor\},S_i\neq y\text{~for all~}0\leq i\leq j-1\Big\}
\]
which is a random set depending on $\{S_0,S_1,\ldots,S_{j-1}\}$. For each $j\geq 1$ such that $\lfloor\log d\rfloor \mid j$,
\[
P\big(S_j=S_{j-1}+e_l\big|S_i,0\leq i\leq j-1\big)=\frac{1}{\lfloor\frac{d}{\log d}\rfloor}
\]
for each $d-\lfloor \frac{d}{\log d}\rfloor+1\leq l\leq d$. For any $x=(x_1,\ldots,x_d)\in \mathbb{Z}^d$, we define
\[
\beta(x)=\sum_{j=d-\lfloor\frac{d}{\log d}\rfloor+1}^{d}x_j,
\]
then it is easy to check that, for each $k\geq 0$ and $k\lfloor\log d\rfloor\leq j\leq (k+1)\lfloor\log d\rfloor-1$,
\[
\beta(S_j)=k.
\]
For each $j$ such that $\lfloor\log d\rfloor\nmid j$, we claim that
\begin{equation}\label{equ 3.4}
|H(j)|\geq 2(d-\lfloor\frac{d}{\log d}\rfloor)-\lfloor\log d\rfloor.
\end{equation}
This is because $\beta(y)=\beta(S_{j-1})$ for each $y$ such that $y-S_{j-1}\in \{\pm e_l:1\leq l\leq d-\lfloor\frac{d}{\log d}\rfloor\}$ while ${\rm card}\{u:\beta(S_u)=\beta(S_j)\}=\lfloor\log d\rfloor$.

For each $k\geq 1$, we use $\vec{S}_k$ to denote the path $(S_0,S_1,\ldots, S_{k})$ on $\mathbb{Z}^d$, then it is easy to check that
\[
\vec{S}_{k\lfloor\log d\rfloor}\in \Gamma_k
\]
for each $k\geq 1$. We denote by $\{V_n\}_{n=0}^{+\infty}$ an independent copy of $\{S_n\}_{n=0}^{+\infty}$ with $V_0=0$, and use $\vec{V}_k$ to denote the path
$(V_0,V_1,\ldots,V_k)$ for each $k\geq 1$, then we define
\[
D(\vec{S},\vec{V})=\bigcup_{k\geq 1}D(\vec{S}_k\lfloor\log d\rfloor,\vec{V}_k\lfloor\log d\rfloor)=\Big\{i: V_i=S_j\text{~for some~}j\geq 0\Big\}
\]
and
\[
F(\vec{S},\vec{V})=\bigcup_{k\geq 1}F(\vec{S}_k\lfloor\log d\rfloor,\vec{V}_k\lfloor\log d\rfloor)=\Big\{i: V_i=S_j\text{~and~}V_{i+1}=S_{j+1}\text{~for some~}j\geq 0\Big\},
\]
where we use $\vec{S}$ and $\vec{V}$ to denote the entire paths of $\{S_n\}_{n\geq 0}$ and $\{V_n\}_{n\geq 0}$ respectively. Here we claim that $|D(\vec{S},\vec{V})|<+\infty$ almost surely under our assumption that $d\geq 20$. The reason is as follows. When $d\geq 20$, $\lfloor \frac{d}{\log d}\rfloor\geq 4$. The path of the latter $\lfloor \frac{d}{\log d}\rfloor$ coordinates of $\vec{S}$ is a $\lfloor\log d\rfloor$ times slower oriented random walk on $\mathbb{Z}^{\lfloor\frac{d}{\log d}\rfloor}$ while former reference shows that two independent oriented simple random walks on $\mathbb{Z}^u$ with $u\geq 4$ collides with each other finitely many times almost surely.

We use $\widetilde{P}$ to denote the probability measure of $\{S_n\}_{n\geq 0}^{+\infty}$ and $\{V_n\}_{n\geq 0}^{+\infty}$ while denote by $\widetilde{E}$ the expectation operator with respect to $\widetilde{P}$, then the following lemma is crucial for us to prove Equation \eqref{equ 2.3}.

\begin{lemma}\label{lemma 3.2}
If $\lambda$ satisfies that
\[
\widetilde{E}\Big(\big(\frac{\lambda+1}{\lambda p}\big)^{|F(\vec{S},\vec{V})|}\big(\frac{2}{p}\big)^{|D(\vec{S},\vec{V})\setminus F(\vec{S},\vec{V})|}\Big)<+\infty,
\]
then $\lambda\geq \lambda_c(d)$.
\end{lemma}
The following lemma is utilized in the proof of Lemma \ref{lemma 3.2}.
\begin{lemma}\label{lemma 3.3}
If $C_1, C_2,\ldots, C_n$ are $n$ arbitrary random events defined under the same probability space such that $P(C_i)>0$ for $1\leq i\leq n$ and $q_1,q_2,\ldots,q_n$ are $n$ positive constants such that
$\sum_{j=1}^nq_j=1$, then
\[
P(\bigcup_{j=1}^{+\infty}C_j)\geq \frac{1}{\sum\limits_{i=1}^n\sum\limits_{j=1}^nq_iq_j\frac{P(C_i\bigcap C_j)}{P(C_i)P(C_j)}}.
\]

\end{lemma}

\proof[Proof of Lemma \ref{lemma 3.3}]

For each $1\leq i\leq n$, we define
\[
Y_i=
\begin{cases}
&\frac{q_i}{P(C_i)} \text{\quad on~} C_i,\\
& 0\text{\quad on~} C_i^c,
\end{cases}
\]
where $C_i^c$ is the complement set of $C_i$, then
\[
\big(E(\sum_{i=1}^nY_i)\big)^2=\big(\sum_{i=1}^n\frac{q_i}{P(C_i)}P(C_i)\big)^2=1^2=1
\]
and
\[
E\Big(\big(\sum_{i=1}^nY_i\big)^2\Big)=\sum_{i=1}^n\sum_{j=1}^nE(Y_iY_i)=\sum\limits_{i=1}^n\sum\limits_{j=1}^nq_iq_j\frac{P(C_i\bigcap C_j)}{P(C_i)P(C_j)}.
\]
According to H\"{o}lder's inequality,
\[
P(\bigcup_{i=1}^nC_i)=P(\sum_{i=1}^nY_i>0)\geq \frac{\big(E(\sum_{i=1}^nY_i)\big)^2}{E\Big(\big(\sum_{i=1}^nY_i\big)^2\Big)}=\frac{1}{\sum\limits_{i=1}^n\sum\limits_{j=1}^nq_iq_j\frac{P(C_i\bigcap C_j)}{P(C_i)P(C_j)}}
\]
and the proof is complete.

\qed

Now we give the proof of Lemma \ref{lemma 3.2}.

\proof[Proof of Lemma \ref{lemma 3.2}]

On the event $\bigcap\limits_{k=1}^{+\infty}\bigcup\limits_{\vec{l}\in \Gamma_k}A_{\vec{l}}$, there exist vertices which have ever been infected with arbitrary large norms and hence $I_t\neq \emptyset$ for any $t>0$. Therefore, according to the Dominated Convergence Theorem,
\begin{equation}\label{equ 3.5}
P_\lambda(I_t\neq \emptyset,\forall~t>0)\geq P_\lambda\big(\bigcap_{k=1}^{+\infty}\bigcup_{\vec{l}\in \Gamma_k}A_{\vec{l}}\big)=\lim_{k\rightarrow+\infty}P_\lambda\big(\bigcup_{\vec{l}\in \Gamma_k}A_{\vec{l}}\big).
\end{equation}
For each $\vec{l}\in \Gamma_k$, we define $q_{\vec{l}}=P_\lambda\big(\vec{S}_{k\lfloor\log d\rfloor}=\vec{l}\big)$, then $\sum_{\vec{l}\in \Gamma_k}q_{\vec{l}}=1$ since $\vec{S}_{k\lfloor\log d\rfloor}\in \Gamma_k$ almost surely. Then by Lemma \ref{lemma 3.3},
\[
P_\lambda\big(\bigcup_{\vec{l}\in \Gamma_k}A_{\vec{l}}\big)\geq \frac{1}{\sum\limits_{\vec{l}\in \Gamma_k}\sum\limits_{\vec{s}\in \Gamma_k}q_{\vec{l}}q_{\vec{s}}\frac{P(A_{\vec{l}}\bigcap A_{\vec{s}})}{P(A_{\vec{l}})P(A_{\vec{s}})}}.
\]
By Equation \eqref{equ 3.3} and Lemma \ref{lemma 3.1},
\[
\frac{P(A_{\vec{l}}\bigcap A_{\vec{s}})}{P(A_{\vec{l}})P(A_{\vec{s}})}\leq \big(\frac{\lambda+1}{\lambda p}\big)^{|F(\vec{l},\vec{s})|}\big(\frac{2}{p}\big)^{|D(\vec{l},\vec{s})\setminus F(\vec{l},\vec{s})|}.
\]
As a result,
\begin{align}\label{equ 3.6}
P_\lambda\big(\bigcup_{\vec{l}\in \Gamma_k}A_{\vec{l}}\big)&\geq \frac{1}{\sum\limits_{\vec{l}\in \Gamma_k}\sum\limits_{\vec{s}\in \Gamma_k}q_{\vec{l}}q_{\vec{s}}\big(\frac{\lambda+1}{\lambda p}\big)^{|F(\vec{l},\vec{s})|}\big(\frac{2}{p}\big)^{|D(\vec{l},\vec{s})\setminus F(\vec{l},\vec{s})|}}\\
&=\frac{1}{\widetilde{E}\Big(\big(\frac{\lambda+1}{\lambda p}\big)^{|F(\vec{S}_{k\lfloor\log d\rfloor},\vec{V}_{k\lfloor\log d\rfloor})|}\big(\frac{2}{p}\big)^{|D(\vec{S}_{k\lfloor\log d\rfloor},\vec{V}_{k\lfloor\log d\rfloor})\setminus F(\vec{S}_{k\lfloor\log d\rfloor},\vec{V}_{k\lfloor\log d\rfloor})|}\Big)},\notag
\end{align}
since $\widetilde{P}\big(\vec{S}_{k\lfloor\log d\rfloor}=\vec{l},~\vec{V}_{k\lfloor\log d\rfloor}=\vec{s}\big)=q_{\vec{l}}q_{\vec{s}}$. According to the Dominated Convergence Theorem,
\begin{align*}
\lim_{k\rightarrow+\infty}&\widetilde{E}\Big(\big(\frac{\lambda+1}{\lambda p}\big)^{|F(\vec{S}_{k\lfloor\log d\rfloor},\vec{V}_{k\lfloor\log d\rfloor})|}\big(\frac{2}{p}\big)^{|D(\vec{S}_{k\lfloor\log d\rfloor},\vec{V}_{k\lfloor\log d\rfloor})\setminus F(\vec{S}_{k\lfloor\log d\rfloor},\vec{V}_{k\lfloor\log d\rfloor})|}\Big)\\
&=\widetilde{E}\Big(\big(\frac{\lambda+1}{\lambda p}\big)^{|F(\vec{S},\vec{V})|}\big(\frac{2}{p}\big)^{|D(\vec{S},\vec{V})\setminus F(\vec{S},\vec{V})|}\Big).
\end{align*}
Then by Equations \eqref{equ 3.5} and \eqref{equ 3.6},
\[
P_\lambda(I_t\neq \emptyset,\forall~t>0)\geq\frac{1}{\widetilde{E}\Big(\big(\frac{\lambda+1}{\lambda p}\big)^{|F(\vec{S},\vec{V})|}\big(\frac{2}{p}\big)^{|D(\vec{S},\vec{V})\setminus F(\vec{S},\vec{V})|}\Big)}>0
\]
and $\lambda\geq \lambda_c(d)$ consequently when
\[
\widetilde{E}\Big(\big(\frac{\lambda+1}{\lambda p}\big)^{|F(\vec{S},\vec{V})|}\big(\frac{2}{p}\big)^{|D(\vec{S},\vec{V})\setminus F(\vec{S},\vec{V})|}\Big)<+\infty.
\]

\qed

We do not check whether $\vec{S}_{k\lfloor\log d\rfloor}$ is uniformly distributed on $\Gamma_k$ since our proof of Lemma \ref{lemma 3.2} does not require this property to hold. We leave this to readers good at calculation.

According to Lemma \ref{lemma 3.2}, we need to give upper bound of $\widetilde{E}\Big(\big(\frac{\lambda+1}{\lambda p}\big)^{|F(\vec{S},\vec{V})|}\big(\frac{2}{p}\big)^{|D(\vec{S},\vec{V})\setminus F(\vec{S},\vec{V})|}\Big)$. We have the following related lemma.

\begin{lemma}\label{lemma 3.4}
There exists a constant $M_1<+\infty$ which does not depend on $d$ such that for any $\theta,\psi>0$,
\[
\widetilde{E}\Big(\theta^{|F(\vec{S},\vec{V})|}\cdot
\psi^{|D(\vec{S},\vec{V})\setminus F(\vec{S},\vec{V})|}\Big)\leq
\psi\sum_{k=0}^{+\infty}\sum_{j=1}^3\Phi_{\theta,\psi}^k(1,j),
\]
where $\Phi_{\theta,\psi}$ is a $3\times 3$ matrix such that
\[
\Phi_{\theta,\psi}=
\begin{pmatrix}
\frac{\big({\lfloor\log d\rfloor}^{\frac{3}{\lfloor\log d\rfloor-1}}\big)\theta}{2(d-\lfloor\frac{d}{\log d}\rfloor)-\lfloor\log d\rfloor}&\frac{\theta}{\lfloor\frac{d}{\log d}\rfloor\lfloor\log d\rfloor^3}&\frac{M_1(\log d)^5\psi}{d}\\
\frac{\big({\lfloor\log d\rfloor}^{\frac{3}{\lfloor\log d\rfloor-1}}\big)\theta}{2(d-\lfloor\frac{d}{\log d}\rfloor)-\lfloor\log d\rfloor}&\frac{\theta}{\lfloor\frac{d}{\log d}\rfloor\lfloor\log d\rfloor^3}&\frac{M_1(\log d)^5\psi}{d}\\
\frac{\big({\lfloor\log d\rfloor}^{\frac{3}{\lfloor\log
d\rfloor-1}}\big)\theta}{2(d-\lfloor\frac{d}{\log
d}\rfloor)-\lfloor\log d\rfloor}&\frac{\theta}{\lfloor\frac{d}{\log
d}\rfloor\lfloor\log d\rfloor^3}&\frac{M_1(\log d)^5\psi}{d}
\end{pmatrix}.
\]
\end{lemma}
We give the proof of Lemma \ref{lemma 3.4} at the end of this section. Now we show how to utilize Lemma \ref{lemma 3.4} to prove Equation \eqref{equ 2.3}.

\proof[Proof of Equation \eqref{equ 2.3}]

For given $r>1$, let $\lambda=\frac{r}{2dp}$, $\theta=\frac{\lambda+1}{\lambda p}$ and $\psi=\frac{2}{p}$, then
\[
\max\Big\{\sum_{j=1}^3\Phi_{\theta,\psi}(i,j):1\leq i\leq 3\Big\}<1
\]
for sufficiently large $d$ according to the definition of $\Phi_{\theta,\psi}$. As a result, by Lemma \ref{lemma 3.4},
\[
\widetilde{E}\Big(\theta^{|F(\vec{S},\vec{V})|}\cdot \psi^{|D(\vec{S},\vec{V})\setminus F(\vec{S},\vec{V})|}\Big)<+\infty
\]
for sufficiently large $d$, where $\theta=\frac{\lambda+1}{\lambda p}$, $\psi=\frac{2}{p}$ and $\lambda=\frac{r}{2dp}$ for $r>1$. Then by Lemma \ref{lemma 3.2},
\[
\lambda_c(d)\leq \lambda=\frac{r}{2dp}
\]
for sufficiently large $d$ and $r>1$. Therefore,
\[
\limsup_{d\rightarrow+\infty}d\lambda_c(d)\leq \frac{r}{2p}
\]
for any $r>1$. Let $r\rightarrow 1$, then the proof is complete.

\qed

At last, we only need to prove Lemma \ref{lemma 3.4}. For this purpose, we introduce some notations and definitions. We use $\tau$ to denote $|D(\vec{S},\vec{V})|$, which is the number of vertices both $\vec{S}$ and $\vec{V}$ visit. For $1\leq i\leq \tau$, we define $t(1)=0$ and
\[
t(i)=\inf\big\{j:j>t(i-1)\text{~and~}j\in D(\vec{S},\vec{V})\big\}.
\]
Note that $t(1)=0$ because $S_0=V_0=O$.  
We divide $\{t(i)\}_{1\leq i\leq \tau}$ into three different types. If $t(i)\in F(\vec{S},\vec{V})$ and $\lfloor\log d\rfloor\nmid \big(t(i)+1\big)$, then we say that $t(i)$ is with type $1$. If $t(i)\in F(\vec{S},\vec{V})$ and $\lfloor \log d\rfloor\mid \big(t(i)+1\big)$, then we say that $t(i)$ is with type $2$. If $t(i)\in D(\vec{S},\vec{V})\setminus F(\vec{S},\vec{V})$, then we say that $t(i)$ is with type $3$. From now on, we assume that $\theta$ and $\psi$ are given. For $i=1,2,3$, we define
\[
\alpha(i)=
\begin{cases}
\theta &\text{~if~}i=1,2,\\
\psi & \text{~if~}i=3.
\end{cases}
\]
For $k\geq 2$ and $1\leq i,j\leq 3$, we use $\Upsilon(k,i,j)$ to denote
\begin{align*}
\alpha(j)\widetilde{P}\Bigg(\tau\geq k+1, t(k)\text{~is with type~}j\Bigg|\vec{S};~t(u),u\leq k; ~V_u,u\leq t(k)\Bigg)
\end{align*}
on the event $\big\{\tau\geq k;~t(k-1)\text{~is with type~}i\big\}$.
For $1\leq j\leq 3$, we define
\[
\nu(j)=\alpha(j)\widetilde{P}\Big(\tau\geq 2, t(1)\text{~is with type~}j\Big|\vec{S}\Big).
\]
For $k\geq 2$ and $1\leq i\leq 3$, we use $b(k,i)$ to denote
\begin{align*}
\psi\widetilde{P}\Bigg(\tau=k\Bigg|\vec{S};~t(u),u\leq k; ~V_u,u\leq t(k)\Bigg)
\end{align*}
on the event $\big\{\tau\geq k;~t(k-1)\text{~is with type~}i\big\}$.
Note that $t(k)$ is with type $3$ when $\tau=k$. For each $m\geq 1$,
we define
\[
W_m=\Big\{\vec{i}=(i_1,i_2,\ldots,i_m)\in \mathbb{Z}^m:i_l\in \{1,2,3\}\text{~for~}1\leq l\leq m\Big\}.
\]
According to the Total Probability Theorem, for $m\geq 2$ and
$\vec{i}=(i_1,\ldots,i_m)\in W_m$,
\begin{align*}
&\widetilde{E}\Big(\theta^{|F(\vec{S},\vec{V})|}\cdot
\psi^{|D(\vec{S},\vec{V})\setminus F(\vec{S},\vec{V})|};
~\tau=m+1,~t(l)\text{~is with type~}i_l\text{~for~}1\leq l\leq
m\Big)\\
&=\widetilde{E}\Big(\nu(i_1)\big[\prod_{k=2}^m\Upsilon(k,i_{k-1},i_k)\big]b(m+1,i_m)\Big)
\end{align*}
and hence
\begin{align}\label{equ 3.7}
&\widetilde{E}\Big(\theta^{|F(\vec{S},\vec{V})|}\cdot \psi^{|D(\vec{S},\vec{V})\setminus F(\vec{S},\vec{V})|}\Big)\notag\\
&=\psi\widetilde{P}(\tau=1)+\widetilde{E}\Big(\sum_{i=1}^3\nu(i)b(2,i)
+\sum_{m=2}^{+\infty}\sum_{\vec{i}\in
W_m}\nu(i_1)\big[\prod_{k=2}^m\Upsilon(k,i_{k-1},i_k)\big]b(m+1,i_m)\Big)\notag\\
&\leq \psi+\psi\widetilde{E}\Big(\sum_{i=1}^3\nu(i)
+\sum_{m=2}^{+\infty}\sum_{\vec{i}\in
W_m}\nu(i_1)\big[\prod_{k=2}^m\Upsilon(k,i_{k-1},i_k)\big]\Big)
\end{align}
since $b(k,i)\leq \psi$. Note that here we utilize the assumption
that $d\geq 20$, which ensures that $\tau<+\infty$ almost surely.

To prove Lemma \ref{lemma 3.4}, we need the following lemma.
\begin{lemma}\label{lemma 3.5}
For each $k\geq 2$ and each $1\leq i\leq 3$,
\begin{equation}\label{equ 3.8}
\Upsilon(k,i,1)\leq \frac{\theta}{2(d-\lfloor\frac{d}{\log d}\rfloor)-\lfloor\log d\rfloor}\text{\quad and \quad}\nu(1)\leq \frac{\theta}{2(d-\lfloor\frac{d}{\log d}\rfloor)-\lfloor\log d\rfloor}.
\end{equation}
For each  $k\geq 2$ and each $1\leq i\leq 3$,
\begin{equation}\label{equ 3.9}
\Upsilon(k,i,2)\leq \frac{\theta}{\lfloor\frac{d}{\log d}\rfloor}\text{\quad and \quad}\nu(2)\leq  \frac{\theta}{\lfloor\frac{d}{\log d}\rfloor}.
\end{equation}
For each  $k\geq 2$ and each $1\leq i\leq 3$, there exists a constant $M_1$ which does not depend on $d,\theta,\psi$ such that
\begin{equation}\label{equ 3.10}
\Upsilon(k,i,3)\leq \frac{\psi M_1(\log d)^2}{d}\text{\quad and \quad}\nu(3)\leq  \frac{\psi M_1(\log d)^2}{d^2}.
\end{equation}
\end{lemma}

\proof

For Equation \eqref{equ 3.8}, conditioned on $\tau\geq k$, there exists a unique $j$ such that $S_j=V_{t(k)}$. Note that $j$ is unique according to the fact that $\vec{S}$ is self-avoiding. Then, $t(k)$ is with type $1$ when and only when $V_{t(k)+1}=S_{j+1}$ and $\lfloor\log d\rfloor\nmid\big(t(k)+1\big)$. As a result,
\[
\Upsilon(k,i,1)\leq \theta\max\Big\{\widetilde{P}\big(V_l=y\big):~y\in \mathbb{Z}^d, \lfloor\log d\rfloor\nmid l\Big\}\leq \frac{\theta}{2(d-\lfloor\frac{d}{\log d}\rfloor)-\lfloor\log d\rfloor}
\]
by Equation \eqref{equ 3.4}. $\nu(1)\leq \frac{\theta}{2(d-\lfloor\frac{d}{\log d}\rfloor)-\lfloor\log d\rfloor}$ follows from a similar analysis.

For Equation \eqref{equ 3.9}, according to a similar analysis with that in the proof of Equation \eqref{equ 3.8},
\[
\Upsilon(k,i,2)\leq \theta\max\Big\{\widetilde{P}\big(V_l=y\big):~y\in \mathbb{Z}^d, \lfloor\log d\rfloor\mid l\Big\}=\frac{\theta}{\lfloor\frac{d}{\log d}\rfloor}.
\]
$\nu(2)\leq  \frac{\theta}{\lfloor\frac{d}{\log d}\rfloor}$ follows from a similar analysis.

For Equation \eqref{equ 3.10},
\[
\Upsilon(k,i,3)\leq \psi \widetilde{P}\Bigg(\exists~j>t(k), j\in D(\vec{S},\vec{V})\Bigg|\vec{S};~t(u),u\leq k; ~V_u,u\leq t(k)\Bigg)
\]
while
\[
\widetilde{P}\Bigg(\exists~j>t(k), j\in D(\vec{S},\vec{V})\Bigg|\vec{S};~t(u),u\leq k; ~V_u,u\leq t(k)\Bigg)=A+B,
\]
where
\[
A=\widetilde{P}\Bigg(\exists~j>t(k), \beta(V_j)=\beta(V_{t(k)})\text{~and~} j\in D(\vec{S},\vec{V})\Bigg|\vec{S};~t(u),u\leq k; ~V_u,u\leq t(k)\Bigg)
\]
and
\[
B=\widetilde{P}\Bigg(\exists~j>t(k), \beta(V_j)>\beta(V_{t(k)})\text{~and~}j\in D(\vec{S},\vec{V})\Bigg|\vec{S};~t(u),u\leq k; ~V_u,u\leq t(k)\Bigg).
\]
If $j\in D(\vec{S},\vec{V})$ for some $j$ such that $j>t(k)$ and $\beta(V_j)=\beta(V_{t(k)})$, then there exists $i$ such that $\beta(S_i)=\beta(V_{{t(k)}})$ and $V_j=S_i$. The value of the function $\beta$ increases by one every $\lfloor\log d\rfloor$ steps of the random walk, hence the number of possible choices of such $(j,i)$ is at most $\lfloor \log d\rfloor^2$. As a result,
\[
A\leq \lfloor \log d\rfloor^2 \max\Big\{\widetilde{P}\big(V_l=y\big):~y\in \mathbb{Z}^d, \lfloor\log d\rfloor\nmid l\Big\}\leq  \frac{\lfloor\log d\rfloor^2}{2(d-\lfloor\frac{d}{\log d}\rfloor)-\lfloor\log d\rfloor}.
\]
Note that in the above equation we utilize the fact that $\lfloor\log d\rfloor\nmid j$ when $j>t(k)$ and $\beta(V_j)=\beta(V_{t(k)})$.

Now we deal with $B$. For each $x=(x_1,\ldots,x_d)\in \mathbb{Z}^d$, we define
\[
\xi(x)=(x_{d-\lfloor\frac{d}{\log d}\rfloor+1},\ldots,x_d)\in \mathbb{Z}^{\lfloor\frac{d}{\log d}\rfloor}.
\]
Then $\{\xi(S_{k\lfloor\log d\rfloor})\}_{k\geq 0}$ and
$\{\xi(V_{k\lfloor\log d\rfloor})\}_{k\geq 0}$ are two independent
oriented simple random walks on $\mathbb{Z}^{\lfloor\frac{d}{\log
d}\rfloor}$ starting at the origin according to our definition of
$\vec{S}$ and $\vec{V}$. It is shown in \cite{Cox1983} that there
exists $M_2>0$ such that two independent oriented simple random walk
on $\mathbb{Z}^d$, both starting at $O$, collide with each other at
least once after leaving $O$ with probability at most $1/d+M_2/d^2$,
where $M_2$ does not depend on $d$. Let
$c=\lfloor\frac{t(k)}{\lfloor\log d\rfloor}\rfloor$, then
\begin{align*}
B&\leq \widetilde{P}\Big(\exists~f>c, \xi(V_{f\lfloor\log
d\rfloor})=\xi(S_{f\lfloor\log d\rfloor})\Big|\xi(V_{c\lfloor\log
d\rfloor})=\xi(S_{c\lfloor\log d\rfloor})\Big)\\
&\leq \frac{1}{\lfloor\frac{d}{\log
d}\rfloor}+\frac{M_2}{\lfloor\frac{d}{\log d}\rfloor^2}.
\end{align*}
As a result,
\begin{align*}
\Upsilon(k,i,3)&= \psi(A+B)\\
&\leq \psi\big(\frac{\lfloor\log
d\rfloor^2}{2(d-\lfloor\frac{d}{\log d}\rfloor)-\lfloor\log
d\rfloor}+\frac{1}{\lfloor\frac{d}{\log
d}\rfloor}+\frac{M_2}{\lfloor\frac{d}{\log d}\rfloor^2}\big)\leq
\frac{\psi M_1(\log d)^2}{d},
\end{align*}
where we can choose $M_1$ which does not depend on $d$. $\nu(3)\leq
\frac{\psi M_1(\log d)^2}{d}$ follows from a similar analysis.

\qed

At last we give the proof of Lemma \ref{lemma 3.4}.

\proof[Proof of Lemma \ref{lemma 3.4}]

According to Equation \eqref{equ 3.7} and Lemma \ref{lemma 3.5},
\begin{align*}
&\widetilde{E}\Big(\theta^{|F(\vec{S},\vec{V})|}\cdot
\psi^{|D(\vec{S},\vec{V})\setminus F(\vec{S},\vec{V})|};~\tau=m+1,
~t(l)\text{~is with type~}i_l\text{~for~}1\leq l\leq m\Big)\\
&=\widetilde{E}\Big(\nu(i_1)\big[\prod_{l=2}^m\Upsilon(l,i_{l-1},i_l)\big]b(m+1,i_m)\Big)\leq
\psi \Lambda(1,i_1)\prod_{l=2}^{m}\Lambda(i_{l-1},i_l)
\end{align*}
for any $\vec{i}=(i_1,i_2,\ldots,i_m)\in W_m$ with $m\geq 2$ and
hence
\begin{align}\label{equ 3.11}
\widetilde{E}\Big(\theta^{|F(\vec{S},\vec{V})|}\cdot
\psi^{|D(\vec{S},\vec{V})\setminus F(\vec{S},\vec{V})|}\Big)&\leq \psi+\psi\sum_{i=1}^3\Lambda(1,i)+\psi\sum_{m=2}^{+\infty}\sum_{\vec{i}\in W_{m}}
\Lambda(1,i_1)\prod_{l=2}^m\Lambda(i_{l-1},i_l)\notag\\
&=\psi\sum_{k=0}^{+\infty}\sum_{j=1}^3\Lambda^k(1,j),
\end{align}
where $\Lambda$ is a $3\times 3$ matrix such that
\[
\Lambda=
\begin{pmatrix}
\frac{\theta}{2(d-\lfloor\frac{d}{\log d}\rfloor)-\lfloor\log d\rfloor}&\frac{\theta}{\lfloor\frac{d}{\log d}\rfloor}&\frac{\psi M_1(\log d)^2}{d}\\
\frac{\theta}{2(d-\lfloor\frac{d}{\log d}\rfloor)-\lfloor\log d\rfloor}&\frac{\theta}{\lfloor\frac{d}{\log d}\rfloor}&\frac{\psi M_1(\log d)^2}{d}\\
\frac{\theta}{2(d-\lfloor\frac{d}{\log d}\rfloor)-\lfloor\log
d\rfloor}&\frac{\theta}{\lfloor\frac{d}{\log d}\rfloor}&\frac{\psi
M_1(\log d)^2}{d}
\end{pmatrix}.
\]
Note that the expression of $\Lambda$ is different with that of
$\Phi_{\theta,\psi}$. We can replace $\Lambda$ by
$\Phi_{\theta,\psi}$ in Equation \eqref{equ 3.11} according to the
following analysis. $\Phi_{\theta,\psi}$ is generated from $\Lambda$
through multiplying the first column of $\Lambda$ by
$\big({\lfloor\log d\rfloor}^{\frac{3}{\lfloor\log
d\rfloor-1}}\big)$, multiplying the second column of $\Lambda$ by
$\frac{1}{{\lfloor\log d\rfloor}^3}$ and multiplying the third
column of $\Lambda$ by ${\lfloor\log d\rfloor}^3$. Between two
adjacent type $2$ moments, there is either at least one type $3$
moment or $\lfloor \log d\rfloor-1$ consecutive type $1$ moments.
Therefore,
\begin{align}\label{equ 3.12}
\Lambda(1,i_1)\prod_{l=2}^{m}\Lambda(i_{l-1},i_l)\leq
\Phi_{\theta,\psi}(1,i_1)\prod_{l=2}^{m}\Phi_{\theta,\psi}(i_{l-1},i_l)
\end{align}
for any $\vec{i}=(i_1,i_2,\ldots,i_m)\in W_m$ with $m\geq 2$ such
that
\begin{equation*}
\widetilde{P}\big(t(l)\text{~is with type~}i_l\text{~for~}1\leq
l\leq m\big)>0.
\end{equation*}
In other words, any $\vec{i}\in W_m$ not satisfying Equation
\eqref{equ 3.12} can not be the vector indicating the types of
$(t(1),\ldots,t(m))$. As a result, we can replace $\Lambda$ by
$\Phi_{\theta,\psi}$ in Equation \eqref{equ 3.11} and the proof is
complete.

\qed

\section{Proof of Equation \eqref{equ 2.4}}\label{section four}
In this section we give the proof of Equation \eqref{equ 2.4}.

\proof[Proof of Equation \eqref{equ 2.4}]

We define
\[
N={\rm card}\big\{x\in \mathbb{Z}^d:~x\in I_t \text{~for
some~}t>0\big\},
\]
then
\[
E_\lambda N=1+\sum_{x\neq O}P_{\lambda}\big(x\in I_t\text{~for
some~}t>0\big)
\]
since $O\in I_0$, where $E_\lambda$ is the expectation operator with
respect to $P_\lambda$. By Equations \eqref{equ 3.2} and \eqref{equ
3.3},
\[
P_{\lambda}\big(x\in I_t\text{~for some~}t>0\big)\leq
\sum_{K=1}^{+\infty}\frac{\lambda^Kp^K{\rm
card}\big\{\vec{l}=(l_0,\ldots,l_K)\in
L_K:~l_K=x\big\}}{(\lambda+1)^K}
\]
for $x\neq O$ and hence
\begin{equation}\label{equ 4.1}
E_\lambda N\leq 1+\sum_{K=1}^{+\infty}\frac{\lambda^Kp^K{\rm
card}\big(L_K\big)}{(\lambda+1)^K}\leq
1+\sum_{K=1}^{+\infty}\frac{\lambda^Kp^K2d(2d-1)^{K-1}}{(\lambda+1)^K},
\end{equation}
since ${\rm card}(L_K)\leq 2d(2d-1)^{K-1}$. By Equation \eqref{equ
4.1}, for sufficiently large $d$ such that $p(2d-1)>1$, $E_\lambda
N<+\infty$ when $\lambda<\frac{1}{(2d-1)p-1}$. Therefore,
$P_\lambda(N<+\infty)=1$ when $\lambda<\frac{1}{(2d-1)p-1}$. On the
event $\{N<+\infty\}$, infective vertices die out when all these $N$
vertices become removed. As a result,
$P_\lambda(I_t=\emptyset\text{~for some~}t>0)=1$ when
$\lambda<\frac{1}{(2d-1)p-1}$ and
\[
\frac{1}{(2d-1)p-1}\leq \lambda_c(d).
\]
Therefore
\[
\liminf_{d\rightarrow+\infty}d\lambda_c(d)\geq
\lim_{d\rightarrow+\infty}\frac{d}{(2d-1)p-1}=\frac{1}{2p}
\]
and the proof is complete.

\qed

\quad

\textbf{Acknowledgments.} The author is grateful to the financial
support from the National Natural Science Foundation of China with
grant number 11501542 and the financial support from Beijing
Jiaotong University with grant number KSRC16006536.

{}
\end{document}